\documentclass[12pt,twoside,reqno]{amsart}

\usepackage{amsmath,amssymb,amsthm}
\usepackage{fancyhdr}
\usepackage{epsfig,color,colordvi}

\usepackage[numbers,sort&compress,square,comma]{natbib}

\newtheorem{thm}{Theorem}[section]
\newtheorem{cor}[thm]{Corollary}
\newtheorem{lem}[thm]{Lemma}
\newtheorem{prop}[thm]{Proposition}

\newtheorem{remarks}[thm]{Remarks}

\theoremstyle{definition}
\newtheorem{defn}{Definition}[section]

\numberwithin{equation}{section} \theoremstyle{remark}

\newcommand{\ee}{\mathbb{E}}

\newcommand{\nn}{\mathbb{N}}
\newcommand{\rr}{\mathbb{R}}
\newcommand{\pp}{\mathbb{P}}

\def\BB{\mathcal B}

\def\FF{\mathcal F}

\def\XX{\mathcal X}

\def\<{\langle}
\def\>{\rangle}

\def\beq{\begin{equation}}
\def\deq{\end{equation}}

\def\bdef{\begin{defn}}
\def\ndef{\end{defn}}

\def\bthm{\begin{thm}}
\def\nthm{\end{thm}}

\def\bprop{\begin{prop}}
\def\nprop{\end{prop}}

\def\brmk{\begin{remarks}}
\def\nrmk{\end{remarks}}

\def\bexa{\begin{exa}}
\def\nexa{\end{exa}}

\def\blem{\begin{lem}}
\def\nlem{\end{lem}}

\def\bcor{\begin{cor}}
\def\ncor{\end{cor}}

\def\bexe{\begin{exe}}
\def\nexe{\end{exe}}

\def\bprf{\begin{proof}}
\def\nprf{\end{proof}}

\def\bdes{\begin{description}}
\def\ndes{\end{description}}

\def\dsp{\displaystyle}

\begin{document}

\author{Guangyu Yang}   
\address{School of
Mathematics and Statistics, Wuhan University, 430072 Hubei, China}
\email{study\_yang@yahoo.com.cn}

\author{Yu Miao}  
\address{School of
Mathematics and Statistics, Wuhan University, 430072 Hubei, China
and School of Mathematics and Information Science, Henan Normal
University, 453007 Henan, China.}
\email{yumiao728@yahoo.com.cn}

\date{September 20, 2006; revised November 13, 2006}

\keywords{Invariance principle; the law of the iterated logarithm;
additive functionals of Markov chains; Dunford-Schwartz operator;
fractional coboundaries; vector-valued martingales.}
\subjclass[2000]{Primary 60F17, secondary 60J10, 60J55}

\begin{abstract}
We prove the Strassen's strong invariance principle for
vector-valued additive functionals of a Markov chain via the
martingale argument and the theory of fractional coboundaries. The
hypothesis is a moment bound on the resolvent.
\end{abstract}

\title[Invariance principle for LIL for Markov chains]
{An invariance principle for the law of the iterated logarithm for
additive functionals of Markov chains}

\maketitle

\section{\textbf{Introduction}} Let $(X_n)_{n\geq0}$ denote a
stationary ergodic  Markov chain defined on a probability space
$(\Omega,\mathcal {F},\pp)$, with values in a measurable space
($\XX,\BB$). Let $Q(x,dy)$ be its transition kernel and $\pi$ the
stationary initial distribution. Furthermore, for $p\geq1$, let
$L^p(\pi)$ denote the space of (equivalence classes of)
$\BB$-measurable functions $g: \XX\to\rr^d$ for some $d\geq1$ and
such that $\|g\|_{p}^p:=\int_\XX |g(x)|^p\pi (dx)<\infty$, and let
$L^p_0(\pi)$ denote the set of $g\in L^p(\pi)$ for which $\int_\XX
gd\pi=0$. Here, $|\cdot|$ denotes the Euclidean norm on $\rr^d$.

 Now, fix $d$ and an $\rr^d$-valued function $g\in L^2_0(\pi)$. For $n\geq0$,
 define
 \beq\label{101}
 S_{n+1}=S_{n+1}(g):=\sum_{i=0}^{n} g(X_i)\;\; {\rm and} \;\;S_0=0.
 \deq
For the question of central limit type results for $S_{n}$, there
have been numerous studies from many angles and under different
assumptions; see Maxwell and Woodroofe \cite{MMMW}, Derriennic and
Lin \cite{DL1} \cite{DL2} and references therein. In this note, we
mainly consider the iterated logarithm type results for $S_n$. Since
the appearance of Strassen's paper \cite{Str}, almost sure
invariance principles for the law of iterated logarithm have been
obtained for a large class of independent and dependent sequence
$(Y_{n})_{n\geq1}$; see Strassen \cite{Str1}, Hall and Heyde
\cite{HH}, and Philipp and Stout \cite{PhS}. Here, the Skorokhod
representation plays an important role.

It is well known that, the law of the iterated logarithm (in short
LIL) is closely related to the central limit theorem (in short CLT)
in some sense. There are several approaches to these
 problems. If the chain is Harris recurrent, then the problems may be
 reduced to the independent case in a certain sense, see Meyn and Tweedie
\cite{MT} and Chen \cite{Ch}. If there exists a solution to
Poisson's equation, $h=g+Qh$,
 then the LIL and CLT problems may be reduced to the martingale case, see also
 Gordin and Lifsic \cite{GL} and Meyn and
Tweedie \cite{MT}. Bhattacharya \cite{Bh} obtained the functional
CLT and LIL for
 ergodic stationary Markov processes by discussing the infinitesimal
 generator.
Wu \cite{Wu} extended the forward-backward martingale decomposition
of Meyer-Zheng-Lyons's type from the symmetric case to the general
stationary situation and gave the Strassen's strong invariance
principle.

Our goal, in this paper, is to consider the problem that $S_n$
satisfies the LIL under some proper conditions. Note that,
Rassoul-Agha and Sepp\"{a}l\"{a}inen \cite{RS} mainly relied on the
invariance principle for vector-valued martingales, so it is likely
to obtain the invariance principle for LIL for the vector-valued
additive functionals of a Markov chain, only if we can develop the
corresponding theory for vector-valued martingales. However, we
encounter the essential difficulties, when considering the
vector-valued martingale, since Monrad and Philipp \cite{MoP} proved
that it is impossible to embed a general $\rr^d$-valued martingale
in an $\rr^d$-valued Gassian process. For the strong approximation
of random sequence taking values in general Banach space, please
refer to Philipp \cite{Phi} and the references given there.

In the present paper, we will take along the lines of Kipnis and
Varadhan \cite{KV} and Maxwell and Woodroofe \cite{MMMW}, to the
case where a solution is not required. Moreover, we identify the
$\lim\sup$ in LIL just the square root of the trace of the diffusion
matrix corresponding to the functional CLT.

Let us explain the organization of this paper. In Section 2, we
state our main results. Section 3 gives the proof of our main
results mentioned in Section 2, which mainly depends on the strong
approximation of vector-valued martingales (see Berger \cite{Be}),
and the theory of fractional coboundaries developed by Derriennic
and Lin \cite{DL}.

\section{\textbf {Main results}}

For a function $h\in L^1(\pi)$, and $\pi$-a.e. $x\in\XX$ define an
operator
\begin{align}
Qh(x)=\int h(y)Q(x,dy).
\end{align}
Obviously, Q is a contraction on $L^p(\pi)$ for $p\geq1$. For
$\varepsilon>0$, let $h_\varepsilon$ be the solution of the equation
$$(1+\varepsilon)h=Qh+g.$$
In fact,
 \beq\label{201}
 h_\varepsilon=\sum_{n=1}^\infty
(1+\varepsilon)^{-n}Q^{n-1}g.
 \deq
 Note that $h_{\varepsilon}\in L^p(\pi)$, if $g\in L^p(\pi)$.
 Let $\pi_1$ be the joint distribution of $X_0$ and $X_1$, so that
 $\pi_1(dx_0, dx_1)=Q(x_0, dx_1)\pi(dx_0)$; denote the
 $L^2$-norm on $L^2(\pi_1)$ by $\|\cdot\|_1$; and let
 $$H_\varepsilon(x_0, x_1)=h_\varepsilon(x_1)-Qh_\varepsilon(x_0)$$
 for $x_0$, $x_1\in\XX$. For any $\varepsilon>0$,
 let
 $$
 M_n(\varepsilon)=\sum_{i=0}^{n-1}H_\varepsilon(X_i,
X_{i+1})
 \;\;{\rm and}\;\;
R_n(\varepsilon)=Qh_\varepsilon(X_0)-Qh_\varepsilon(X_n),
$$
 hence, by the simple computation,
 \beq\label{202}
 S_n(g)=M_n(\varepsilon)+\varepsilon
 S_n(h_\varepsilon)+R_n(\varepsilon).
 \deq
For convenience, we summarize the results of Maxwell and Woodroofe
\cite{MMMW} as the following theorem.

\vskip5pt\noindent \textbf{Theorem MW} {\it Assume that $g\in
L^2_{0}(\pi)$ and that there exists an $\alpha\in(0,1/2)$ such that
\begin{align}\label{Moment}
\|\sum_{i=0}^{n-1}Q^{i}g\|_{2}=O(n^{\alpha}).
\end{align}
Then we have
\begin{enumerate}
 \item The limit
$H=\lim_{\varepsilon\rightarrow0^{+}}H_{\varepsilon}$ exists in
$L^2(\pi_{1})$. Moreover, if one defines
\begin{align*}
M_n =\sum_{i=0}^{n-1}m_{i},
\end{align*}
where $m_{i}=H(X_{i},X_{i+1})$, then $(m_{n})_{n\geq0}$ is a
stationary and ergodic $\pp$-square integrable martingale difference
sequence, with respect to the filtration
$\{\FF_{n}=\sigma(X_{0},\cdots,X_{n})\}_{n\geq0}$;

\item $\|h_{\varepsilon}\|_2=O(\varepsilon^{-\alpha})$, and if
$R_{n}=S_{n}-M_{n}=M_{n}(\varepsilon)-M_{n}+\varepsilon
S_{n}(h_{\varepsilon})+R_{n}(\varepsilon)$, then
\begin{align*}
\ee(|R_{n}|^2)=O(n^{2\alpha}).
\end{align*}
\end{enumerate}}
\begin{remarks}
Furthermore, if there exists a $p>2$ such that $g\in L^{p}(\pi)$,
then there exists a $q\in (2, p)$ such that $H\in L^q(\pi)$ and
$(M_{n})_{n\geq1}$ is an $L^{q}$-martingale; see the Theorem 1 of
Rassoul-Agha and Sepp\"{a}l\"{a}inen \cite{RS}.
\end{remarks}
For introducing our main results, we need give a few more notations.
Let $C([0,1],\rr^d)$ be the Banach space of continuous maps from
$[0,1]$ to $\rr^d$, endowed with the supremum norm $|\|\cdot|\|$,
using the Euclidean norm in $\rr^d$. Denote $K$ the set of
absolutely continuous maps $f\in C([0,1],\,\rr^d)$, such that
\begin{align*}
f(0)=0, \ \ \ \ \int_{0}^{1}|\dot{f}(t)|^2dt\leq1,
\end{align*}
where, $\dot{f}$ denotes the derivative of $f$ determined almost
everywhere with respect to Lebesgue measure. Obviously, $K$ is
relatively compact and closed.

Define
\begin{align*}
\xi_{n}(t)=(2n\log\log n)^{-1/2}[S_{k}+(nt-k)g(X_{k})]
\end{align*}
for $t\in[\frac{k}{n},\frac{k+1}{n}),\;k=0,1,2,\cdots,n-1.$ In order
to avoid difficulties in specification, we adopt the convention that
$\log\log x=1$, if $0<x\leq e^e$. Then, $\xi_{n}$ is a random
element with values in $C([0,1],\,\rr^d)$.

After these preparations, we are now in a position to state our main
results.
\begin{thm}\label{thm1}
Let $g\in L^p_0(\pi)\,(p>2)$ and assume that there exists an
$\alpha\in(0,1/2)$ for which (\ref{Moment}) is satisfied. Then, the
sequence of functions $(\xi_{n}(\cdot), \,n\geq1)$ is relatively
compact in the space $C([0,1],\,\rr^d)$, and the set of its limit
points as $n\rightarrow\infty$, coincides with $\sqrt{{\rm
tr}(\mathfrak{D})}K$, where {\rm tr($\cdot$)} denotes the trace
operator of a matrix and $\mathfrak{D}=\ee(M_1M_1^t)=\int
HH^td\pi_{1}$ is the diffusion matrix corresponding to the
functional central limit theorem; see Rassoul-Agha and
Sepp\"{a}l\"{a}inen \cite{RS}.
\end{thm}
\begin{thm}\label{thm2}
Let $g\in L^p_0(\pi)\,(p>2)$ and assume that there exists an
$\alpha\in(0,1/2)$ for which (\ref{Moment}) is satisfied. Then
\begin{align}
\limsup|S_{n}|/\sqrt{2n\log\log n}=\sqrt{\rm tr(\mathfrak{D})}
\;\;\;\;\;\pp-a.s.
\end{align}
\end{thm}
\begin{remarks}
In fact, ${\rm tr}(\mathfrak{D})=\|H\|_{1}^2$. And particularly, if
putting $d=1$, we can obtain the main results of Miao and Yang
\cite{MY}.
\end{remarks}
\begin{remarks} For $n\geq0$, define $
S_n^*=S_n-\ee_{X_{0}}S_{n}$, since the Theorem 3 of Rassoul-Agha and
Sepp\"{a}l\"{a}inen \cite{RS}, the above Theorem 2.1 and Theorem 2.2
also hold for $S_{n}^*$.
\end{remarks}
\section{\textbf{Proof of main results}}
\subsection{ Proof of Theorem \ref{thm1}}

\begin{proof}
For $0\leq t\leq1$, define
\begin{align*}
\zeta_{n}(t)=(2n\log\log n)^{-1/2}M_{[nt]},\\\eta_{n}(t)=(2n\log\log
n)^{-1/2}B(nt),
\end{align*}
where, $M_{n}$ is as defined in Section 2 and $B(\cdot)$ is an
$\rr^d$-valued Brownian motion with mean 0 and diffusion matrix
$\mathfrak{D}$. Theorem 1 of Strassen \cite{Str} shows that
$(\eta_{n}(\cdot))_{n\geq1}$ is relatively compact and the set of
its limit points coincides with $\sqrt{{\rm tr}(\mathfrak{D})}K$.

Notice that by the part (1) of Theorem MW, $(M_n)_{n\geq1}$ is a
square integrable martingale with strictly stationary increments.
Moreover,
\begin{align}
\ee(<u,m_0>^2)<\infty\;\;{\rm and} \;\;\ee(<u,m_0>)=0, \;\;{\rm
for\;all}\;u\in\rr^d,
\end{align}
where, $<\cdot,\cdot>$ denotes the inner product in $\rr^d$.
Therefore, Corollary 4.1 of Berger \cite{Be} implies that,

{\it Without changing its distribution, one can redefined the
sequence $(M_{n})_{n\geq1}$ on a new probability space
$(\hat{\Omega},\hat{\mathcal {F}},\hat{\pp})$ on which there exists
an $\rr^d$-valued Brownian motion $(B(t))_{t\geq0}$ with mean 0 and
diffusion matrix $\mathfrak{D}$ such that
\begin{align}
|M_{[t]}-B(t)|=o((t\log\log t)^{-1/2}),
~~~~\hat{\pp}-a.s.~~~~(as\;t\rightarrow\infty).
\end{align}
where, $\mathfrak{D}=\lim_{n\to\infty}n^{-1}Cov(M_n)$.}
\begin{remarks}
Birkhoff-Khinchin's ergodic theory and together with the simple
calculation shows that, $\mathfrak{D}=\ee(M_1M_1^t)=\int
HH^{t}d\pi_{1}$, is the diffusion matrix corresponding to the
functional central limit theorem; see also Rassoul-Agha and
Sepp\"{a}l\"{a}inen \cite{RS}.
\end{remarks}

That is to say,
$$
\sup_{0\leq t\leq1}|M_{[nt]}-B(nt)|=o((2n\log\log
n)^{-1/2}),\;\;\;\;\;\pp-a.s.
$$

 Hence,
\begin{align*}
|\|\zeta_{n}-\eta_{n}|\|=&(2n\log\log
n)^{-1/2}\sup_{0\leq t\leq1}|M_{[nt]}-B(nt)|\\
=&o(1),\ \ \ \ \ \ \ \ \ \ \ \ \ \ \ \pp-a.s.
\end{align*}

Define
\begin{align*}
\tilde{\zeta}_{n}(t)=(2n\log\log n)^{-1/2}[M_{k}+(nt-k)m_{k}]
\end{align*}
for $t\in[\frac{k}{n},\frac{k+1}{n}),\;k=0,1,2,\cdots,n-1.$ Then
$\tilde{\zeta}_{n}\in C([0,1],\,\rr^{d})$ and
\begin{align*}
\sup_{t\in[0,1)}|\zeta_{n}(t)-\tilde{\zeta}_{n}(t)|=(2n\log\log
n)^{-1/2}\max_{0\leq k\leq n-1}|m_{k}|.
\end{align*}
Next, we give the order estimation of $\max_{0\leq k\leq
n-1}|m_{k}|$.

\vskip5pt\noindent\textbf{Lemma M}\label{lemM} {\it ( See M\'oricz
\cite{Mor} ) Let $p
> 0$ and $\beta
> 1$ be two positive real numbers and $Z_i$ be a sequence of random
variables. Assume that there are nonnegative constants $a_j$
satisfying
\begin{align}
\ee|\sum_{j=1}^iZ_j|^p\leq (\sum_{j=1}^ia_j)^\beta,
\end{align}
for $1\leq i\leq n$. Then
\begin{align}
\ee(\max_{1\leq i\leq n}|\sum_{j=1}^iZ_j|^p)\leq
C_{p,q}(\sum_{i=1}^na_i)^\beta,
\end{align}
for some positive constant $C_{p,\beta}$ depending only on $p$ and
$\beta$. }

\begin{lem}
For any enough large $n$, there exists a positive constant C such
that
\begin{align}
\ee(\max_{1\leq i\leq n}|m_i|^q)\leq C \ee|m_1|^q.
\end{align}
\end{lem}
\begin{proof}
Since the part (1) of Theorem MW and Remarks 2.1, let $\ee
|m_1|^q=a^2(q)$ and for any $k\geq1$, we have the following
relations,
\begin{align}
\ee|m_k|^q\leq (\sum_{i=1}^ka_i)^2,
\end{align}
where $a_{1}=a(q)$ and $a_{i}=0$ for $2\leq i\leq k$. Hence, by
Lemma M, there exists a constant $C>0$, such that
\begin{align}
\ee(\max_{1\leq i\leq n}|m_i|^q)\leq C(\sum_{i=1}^n
a_i)^2=C\ee(|m_1|^q).
\end{align}

This completes the proof of the lemma.
\end{proof}

For any $\epsilon>0$, Lemma 3.2 immediately yields,
\begin{align*}
\pp(\max_{0\leq k\leq n-1}|m_{k}|\geq\epsilon(2n\log\log
n)^{1/2})=O((n\log\log n)^{-q/2}).
\end{align*}
By Borel-Cantelli's lemma, we have
\begin{align}
(2n\log\log n)^{-1/2}\max_{0\leq k\leq
n-1}|m_{k}|=o(1),\;\;\;\;\;\pp-a.s.
\end{align}
Hence,
\begin{align*}
|\|\zeta_{n}-\tilde\zeta_{n}|\|=\sup_{t\in[0,1)}|\zeta_{n}(t)-\tilde{\zeta}_{n}(t)|=o(1),\;\;\;\;\;\pp-a.s.
\end{align*}

 The above discussions immediately yield the following claim:
\vskip5pt\noindent{\it $(\tilde{\zeta}_{n}(\cdot),n\geq1)$ is
relatively compact and the set of its limit points coincides with
$\sqrt{{\rm tr}(\mathfrak{D})}K$. }

\vskip5pt Now, we turn to deal with the neglectable term $R_{n}$ in
the sense of functional LIL. Firstly, let us recall the concept of
Dunford-Schwartz (DS) operator; see Derriennic and Lin \cite{DL}. We
call $T$ a DS operator on $L^1$ of a probability space, if $T$ is a
contraction of $L^1$ such that $\|Tf\|_\infty\leq \|f\|_\infty$ for
every $f\in L^\infty$. If $\theta$ is a measure preserving
transformation in a probability space $(\Omega, \Sigma, \mu)$, then
the operator $Tf=f\circ\theta$ is a DS operator on $L^1(\mu)$. More
generally, any Markov transition operator $P$ with an invariant
probability measure yields a positive DS operator.

\vskip5pt\noindent\textbf{Lemma DL}\label{lemDL} {\it( See
Derriennic and Lin \cite{DL}) \begin{enumerate}
\item Let $T$ be a contraction in a Banach space
 $X$, and let $0<\beta<1$. If \\$\dsp
 \sup_{n}\|\frac{1}{n^{1-\beta}}\sum_{k=1}^nT^ky\|<\infty$,
 then $y\in (I-T)^\alpha X$ for every $0<\alpha<\beta$.

 \item Let $T$ is a DS operator in $L^1(\mu)$ of a
 probability space, and fix $1<p<\infty$, with dual $q=p/(p-1)$. Let
 $0<\alpha<1$, and $f\in(I-T)^\alpha L^p$. If
 $\alpha>1-\frac{1}{p}=\frac{1}{q}$, then
 $\dsp \frac{1}{n^{1/p}}\sum_{k=0}^{n-1}T^k f\to 0 $ a.e.
\end{enumerate}
}

To apply the above lemma, i.e., the theory of fractional
coboundaries named by Derriennic and Lin \cite{DL}, we need
construct a DS operator. On $\XX\times\XX$, define
$$
f(x_0, x_1)=g(x_0)-H(x_0, x_1),
$$
then
 \begin{align}\label{208}
 R_n=&S_n-M_n\nonumber\\=&\sum_{i=0}^{n-1} [g(X_{i})-H(X_{i},
 X_{i+1})]\nonumber\\
 =&\sum_{i=0}^{n-1}
f(X_{i}, X_{i+1}).
 \end{align}
Let $\theta$ be the shift map on the path space $\XX^\nn$ for the
Markov chain which is a contraction on $L^2(\pp)$. Hence, $\theta$
is a DS operator. For a sequence $x=(x_{i})_{i\in\nn}\in\XX^\nn$,
define $F(x)=f(x_{0},
  x_{1})$, then we have
$$
F\in L^2(\pp)\;\;\;{\rm
and}\;\;\;R_n=\sum_{k=0}^{n-1}F\circ\theta^k.
$$

From the part (2) of Theorem MW, there exists a constant
$1/2<\beta<1-\alpha$, such that
 \beq\label{210}
\sup_{n}
\Big\|\frac{1}{n^{1-\beta}}\sum_{k=0}^{n-1}F\circ\theta^k\Big\|<\infty,
 \deq
Since the part (1) of Lemma DL and $0<\alpha<1/2$, we have
$F\in(I-\theta)^\eta L^2(\pp)$, for some $\eta\in(1/2, 1-\alpha)$.
By the part (2) of Lemma DL, we
 have
 $$\frac{1}{n^{1/2}}R_n\to 0,\ \ \ \ \ \ \ \ \ \pp-a.s.$$
Furthermore, applying an elementary property of real convergent
sequences, we immediately get
$$
\max_{0\leq k\leq n}|R_{k}|=o((2n\log\log n)^{1/2}),
\;\;\;\;\;\pp-a.s.
$$
Consequently,
\begin{align}
(2n\log\log n)^{-1/2}\sup_{0\leq t\leq1}|R_{[nt]}|\to
0,\;\;\;\;\;\pp-a.s.
\end{align}
From above discussions, we complete the proof of Theorem \ref{thm1}.
\end{proof}

\subsection{Proof of Theorem \ref{thm2}}

\begin{proof}
Here, we take along the lines of the proof of Theorem 4.8 in Hall
and Heyde \cite{HH}. Let $\{e_{i}\}_{i=1}^{d}$ the canonical basis
of $\rr^d$. For any $\rr^d$-valued function $f$, denote
$f=(f_{1},f_{2},\cdots,f_{d})^t$. By the definition of $K$, we have,
for any $f\in\sqrt{tr(\mathfrak{D})}K$,
\begin{align}
|f(t)|^2&=\sum_{i=1}^{d}(\int_{0}^{t}\dot{f}_{i}(s)ds)^2\nonumber\\
&\leq\sum_{i=1}^{d}(\int_{0}^{t}\dot{f}_{i}(s)^2ds)\int_{0}^{t}1ds\leq
tr(\mathfrak{D})t
\end{align}
where, the first inequality by the Cauchy-Schwartz's inequality. So,
$|f(t)|\leq \sqrt{tr(\mathfrak{D})t}$. It follows that
$\sup_{t\in[0,1]}|f(t)|\leq\sqrt{tr(\mathfrak{D})}$.
 Hence, by Theorem 2.1,
\begin{align}
\limsup\sup_{t\in[0,1]}|\xi_{n}(t)|\leq\sqrt{tr(\mathfrak{D})},\;\;\;\;\;\pp-a.s.
\end{align}
and setting $t=1$,
\begin{align}
\limsup|S_{n}|/\sqrt{2n\log\log n}\leq
\sqrt{\textrm{tr}(\mathfrak{D})},\;\;\;\;\;\pp-a.s.
\end{align}

On the other hand, we put
$f(t)=t\sqrt{\frac{tr(\mathfrak{D})}{d}}\sum_{i=1}^{d}e_{i}$, $t\in
[0,1]$. Then, $f\in\sqrt{tr(\mathfrak{D})}K$ and so for
$\pp-a.s.\,\omega$, there exists a sequence $n_{k}=n_{k}(\omega)$,
such that
\begin{align}
\xi_{n_{k}}(\cdot)(\omega)\stackrel{|\|\cdot|\|}{\longrightarrow}
f(\cdot).
\end{align}
 Particularly,
$f(1)=\sqrt{\frac{tr(\mathfrak{D})}{d}}\sum_{i=1}^{d}e_{i}$,
$|\xi_{n_{k}}(1)(\omega)|\longrightarrow |f(1)|$. That is to say,
\begin{align}
|S_{n_{k}}|/\sqrt{2{n_{k}}\log\log
n_{k}}=\sqrt{\textrm{tr}(\mathfrak{D})},\;\;\;\;\;\pp-a.s.
\end{align}
This completes the proof of Theorem \ref{thm2}.
\end{proof}

\section*{\small Acknowledgements}
{\small The authors wish to thank Prof. Y. Derriennic and Prof. M.
Lin for sending the key reference paper \cite{DL}.}

\begin{small}

\end{small}
\end{document}